\documentclass[12pt ,twoside ,a4paper]{article}
\usepackage{enumitem}
\setlist{noitemsep, topsep= 0pt, parsep=1pt, partopsep=1pt, leftmargin=0.5cm}
\normalfont
\usepackage{multirow}
\usepackage{natbib}
\usepackage{amsthm}
\newtheorem{theorem}{Theorem}[section]
\newtheorem{lemma}[theorem]{Lemma} 
\newtheorem{remark}[theorem]{Remark}
\newtheorem{definition}[theorem]{Definition}
\usepackage[a4paper, total={185mm,290mm}, bottom= 39mm,  top= 20mm]{geometry}
\usepackage{amsmath}

\usepackage[pagebackref]{hyperref}
\usepackage[thinlines]{easytable}
\usepackage{hyperref}
\hypersetup{colorlinks=true, linkcolor=blue, urlcolor={blue},
   pdfstartview={FitV},unicode,breaklinks=true, citecolor = blue}
\usepackage{amsmath}
\usepackage[font=scriptsize]{caption}
\usepackage{array}

\usepackage[utf8]{inputenc}
\usepackage{amssymb}
\usepackage{breqn}

\title{\textbf{On limiting distribution of U-statistics based on associated random variables}}
\vspace{-3.0pt}
\author{$\text{ Mansi Garg}$ and $\text{Isha Dewan}$ \vspace{-0.1in}\\
 ${}$Indian Statistical Institute \vspace{-0.1in}\\New Delhi-110016 (India)\vspace{-0.1in}\\ mansibirla@gmail.com and ishadewan@gmail.com $\vspace{-0.3in}$ 
 }
\begin{document}
\fontsize{15}{20}
\date{}
\maketitle
\abstract{Let $\{X_n, n \ge 1\}$ be a sequence of stationary associated random variables. We discuss another set of conditions under which a central limit theorem for U-statistics based on $\{X_n, n \ge 1\}$ holds. We look at U-statistics based on differentiable kernels of degree 2 and above. We also discuss some applications.  }\\
\\
\textbf{Keywords}: {\it Associated random variables;  U-statistics; Central limit theorem; Skewness; Kurtosis.}

\numberwithin{equation}{section}\label{section1}
\section{Introduction}
\setlength\parindent{24pt}

In this paper, we obtain  another set of assumptions under which a central limit theorem  for  U-statistics based on stationary associated random variables holds. We look at non-degenerate U-statistics based on differentiable (component-wise monotonic or non-monotonic)  kernels of any finite degree $k \ge 2$.  The proof requires  relatively non-restrictive assumptions. We also illustrate some applications of our results using examples. Apropos our discussion, we give the following.

\begin{definition}  (\cite{MR0217826}) A finite collection of random variables $\{X_j,  1 \le j \le n \}$  is said to be associated, if for any choice of component-wise nondecreasing functions $h$, $g$ $:$ ${\mathbb{R}}^n \rightarrow {\mathbb{R}}$, we have
\begin{equation*}
Cov(h(X_1,\ldots, X_n), g(X_1,\ldots, X_n)) \ge 0,
\end{equation*}
whenever it exists. An infinite collection of random variables $\{X_j,  j \ge 1 \}$ is associated if every finite sub-collection is associated. 
\end{definition}
Any set of independent random variables is associated. Nondecreasing functions
of associated random variables are associated (cf. \cite{MR0217826}).  A detailed presentation of the asymptotic results and examples relating to associated sequences can be found in \cite{Bulinski:1701595},  \cite{oliveira2012asymptotics} and \cite{MR3025761}.

For the rest of the paper, assume that  $\{X_n, n \ge 1 \}$ is a sequence of stationary associated random variables with $F$ as the marginal distribution function of $X_1$. We next briefly discuss  existing results on central limit theorem of U-statistics based on $\{X_n, n \ge 1 \}$.

\cite{dewanprakasarao2001} gave a central limit theorem for degenerate and non-degenerate U-statistics  using an orthogonal expansion of the underlying kernel. \cite{Dewan20029} and its corrigendum \cite{Dewan2015147}  obtained a central limit theorem for U-statistics with differentiable kernels of degree 2, using Hoeffding's decomposition. The limiting distribution of U-statistics based on differentiable kernels can also be obtained using the results of  \cite{beutner2012, beutner2014} and \cite{Garg2015209}. We have discussed the difference in our assumptions with an example in section $\ref{section4}$.
  
The paper is organized as following. In section $\ref{section2}$, we state a few results and definitions which will be required to prove our main results. Limiting distribution of non-degenerate U-statistics based on $\{X_n, n \ge 1 \}$ is given in section $\ref{section3}$. In section $\ref{section4}$, we apply our results to discuss the asymptotic distribution of estimators of second, third and the fourth central moments. We also discuss estimators of skewness and kurtosis, when the underlying sample is from $\{X_n, n \ge 1 \}$.  

\section{Preliminaries} \label{section2}
\setlength\parindent{24pt}
     In this section, we give results and definitions which will be needed to  prove our main results given in section $\ref{section3}$.
     \begin{definition}\label{hoeff} Hoeffding's decomposition for U-statistics based on a symmetric measurable function $\rho : \mathbb{R}^2 \to \mathbb{R}$.    Define the U-statistic, $U_n$, by
\begin{equation}
U_n = {{n}\choose{2}}^{-1}\sum_{1 \le i < j \le n} {\rho}(X_{i}, X_{j}).
\end{equation}
Let $\theta = \int\limits_{\mathbb{R}^2}{\rho}(x_1, x_2)~{dF({x_1})dF({x_2})}$ \text{ and} ${\rho}_{1}({x_1}) = \int\limits_{\mathbb{R}}{\rho}(x_1, x_2)~dF({x_{2}})$. Further, let
 \begin{align}\label{h2defn}
  h^{(1)}({x_1}) =  {\rho}_1({x_1}) - \theta \text{ and }  h^{(2)}({x_1}, {x_2}) =  {\rho}({x_1}, {x_2}) - {\rho}_1({x_1}) - {\rho}_1({x_2}) + \theta.
  \end{align}
   Then, the Hoeffding's decomposition  for $U_n$ is 
   \begin{equation}
   U_n = \theta + 2H_n^{(1)} + H_n^{(2)},
   \end{equation}
    where $H_{n}^{(j)}$ is the U-statistic of degree $j$ based on the kernel $h^{(j)}$, $j = 1, 2$. When ${X_j,  1 \le j \le n}$ are $i.i.d.$, $E(U_n) = \theta$. 
    
    \begin{remark}An extension of the Hoeffding's decomposition for U-statistics of a finite  degree $k > 2$ can be found in \cite{1990u}.
     \end{remark}
    \end{definition}

\begin{lemma}
 (\cite{newman1980}) \label{covnew1} Let $X$ and $Y$ be two associated random variables with $E(X^2)< \infty$ and $E(Y^2) < \infty$. Let $f$ and $g$ be differentiable functions with $\underset{x}{sup}\left|f'(x)\right| < \infty$ and $\underset{y}{sup}\left|g'(y)\right| < \infty$. Then,
\begin{align*}
 Cov(f(X), g(Y)) &= \int_{\mathbb{R}^2}f'(x)g'(y)[P(X \le x, Y \le y) - P(X \le x)P(Y \le y)] dxdy\\
  & \le \underset{x}{sup}\left|f'(x)\right|  \underset{y}{sup}\left|g'(y)\right|Cov(X, Y).
\end{align*}
\end{lemma}

\begin{lemma}\label{lebowitz} (\cite{lebowitz1972})  Define, for $A$ and $B$, subsets of $\lbrace{1, 2, ..., n}\rbrace$ and real $x_j$$'$s,
 \begin{equation*}
H_{A, B}(x_j, j \in A\cup B) = P[ X_j > x_j, j \in A\cup B] - P[ X_k > x_k, k \in A]P[ X_l > x_l, l \in B].
\end{equation*}
If the random variables  ${X_1, X_2, ..., X_n}$ are associated, then
\begin{equation}
0 \le H_{A, B} \le \sum_{i \in A} \sum_{j \in B} H_{\lbrace{i}\rbrace, \lbrace{j}\rbrace}.
\end{equation}
\end{lemma}

\begin{definition} (\cite{MR789244})\label{covnew2}  Let $g$ and $\tilde{g}$ be two real-valued functions on $\mathbb{R}^{m}$, for some $m \in \mathbb{N}$. $g \ll \tilde{g}$ iff $\tilde{g} + g$ and $\tilde{g} - g$ are both coordinate-wise nondecreasing. If $g \ll \tilde{g}$, then $\tilde{g}$ will be a coordinate-wise nondecreasing function. 
\end{definition}

{\begin{lemma} \label{cltfass}(\cite{MR789244}) For each $j$, $j \ge 1$, let $Y_j = f(X_j)$ and $\tilde{Y_j} = \tilde{f}(X_j)$. Suppose that $f \ll \tilde{f}$. Define $\sigma^2 = Var (Y_1) + 2\sum_{j=2}^{\infty}Cov(Y_1, Y_j)$. Let $\sigma^2 > 0$ and $ \sum_{j=1}^{\infty}Cov(\tilde{Y_1},\tilde{Y_j}) < \infty$.
Then,
\begin{equation}
\frac{1}{\sqrt{n}\sigma}\sum_{j=1}^{n}(Y_j - E(Y_j)) \xrightarrow {\mathcal{ L}} N(0, 1) \: \text{as}\: {n \to \infty}.
\end{equation}
\end{lemma}}

In the following, $g \ll_A \tilde{g}$ if $g \ll \tilde{g}$ and both $g$ and $\tilde{g}$ depend only on $x_j's$ with $j \in A$. $A$ is a finite subset of $\lbrace{k, k \ge 1}\rbrace$.

\begin{lemma} (\cite{MR789244}) \label{covnew3} Let $g_1 \ll_A \tilde{g}_1$ and $g_2 \ll_A \tilde{g}_2$. Then,
\begin{equation}
|Cov(g_1(X_1, X_2, ...), g_2(X_1, X_2, ...))| \le Cov(\tilde{g}_1(X_1, X_2, ...), \tilde{g}_2(X_1, X_2, ...)).
\end{equation}
\end{lemma}

\section{Asymptotic distribution of  U-statistics}\label{section3}
\setlength\parindent{24pt}

The main result of this section, Theorem $\ref{theorem1chap1}$, gives a central limit theorem for non-degenerate U-statistics with  differentiable kernels of degree 2 based on $\{X_n, n \ge 1\}$. The extension of this theorem to U-statistics with kernels of a general finite degree $k > 2$ is also discussed. Before proceeding to the proof of the main theorem, we discuss the following lemma.  

Assume all the expectations and derivatives defined in this section exist, and $C$ is a generic positive constant in the sequel.

\begin{lemma}\label{chris} Let the random variables $Z_1, Z_2, Z_3, Z_4$ be identically distributed with $F$ as the marginal distribution function of $Z_1$ and $E[Z_1^2] < \infty$. Define $z(u, v, x, y) =  h^{(2)}(u, v)h^{(2)}(x, y)$, where the kernel $h^{(2)}(. , .)$ is defined in  $(\ref{h2defn})$. Then, for  ${Z_1', Z_2', Z_3', Z_4'}$    $i.i.d.$ random variables with  marginal distribution function $F$ and  independent of $Z_1, Z_2, Z_3, Z_4$, 
 \begin{align} \label{chrisequation}
  &E\Big[ z(Z_1, Z_2, Z_3, Z_4) - z(Z'_1, Z_2, Z_3, Z_4) \Big] \nonumber \\
     & = E\Big[  \int_{\mathbb{R}^4} \Big(I(Z_1 > u) - I(Z'_1 > u) \Big)\Big(I(Z_2 > v) - I(Z'_2 > v) \Big)\Big(I(Z_3 > x) - I(Z'_3 > x) \Big)\nonumber \\ &\hspace{1in}\Big(I(Z_4 > y) - I(Z'_4 > y) \Big)z''''(u, v, x, y)dudvdxdy\Big],
 \end{align}
 where $z''''(u, v, x, y) = \Big(\frac{\partial^4 z(t_1, t_2, t_3, t_4)}{\partial dt_1\partial dt_2\partial dt_3\partial dt_4 }|_{t_1 = u, t_2 = v, t_3 = x, t_4 = y}\Big)$.
\end{lemma}
\begin{proof} From the results discussed in \cite{CHRISTOFIDES2004138} we get
\begin{align}\label{eq1chris1chap1were}
& z(Z_1, Z_2, Z_3, Z_4) - z(Z'_1, Z_2, Z_3, Z_4)  = \int_{\mathbb{R}} \Big(I(Z_1 > u) - I(Z'_1 > u) \Big)z'(u, Z_2, Z_3, Z_4)du, 
\end{align}
where $z'(u, Z_2, Z_3, Z_4) = \Big(\frac{\partial z(t, Z_2, Z_3, Z_4)}{\partial dt}|_{t = u}\Big)$. Let $z''(u, v, Z_3, Z_4) = \Big(\frac{\partial^2 z(t_1, t_2, Z_3, Z_4)}{\partial dt_1\partial dt_2}|_{t_1 = u, t_2 = v}\Big)$. Then, 
\begin{align}
& z(Z_1, Z_2, Z_3, Z_4) - z(Z'_1, Z_2, Z_3, Z_4) \nonumber \\
   & = \int_{\mathbb{R}} \Big(I(Z_1 > u) - I(Z'_1 > u) \Big)\Big(z'(u, Z_2, Z_3, Z_4)- z'(u, Z'_2, Z_3, Z_4)\Big)du \nonumber \\
    &\hspace{1in} +  \int_{\mathbb{R}} \Big(I(Z_1 > u) - I(Z'_1 > u) \Big)z'(u, Z'_2, Z_3, Z_4)du\nonumber \\
    & = \int_{\mathbb{R}^2} \Big(I(Z_1 > u) - I(Z'_1 > u) \Big)\Big(I(Z_2 > v) - I(Z'_2 > v) \Big)z''(u, v, Z_3, Z_4)dvdu \nonumber \\
 &\hspace{1in} +  z(Z_1, Z'_2, Z_3, Z_4) - z(Z'_1, Z'_2, Z_3, Z_4).\label{obxchap3}
 \end{align}
 The last equality follows by repeating the argument of $(\ref{eq1chris1chap1were})$. 
  Taking expectations, we get 
 \begin{align} \label{eq1chap1chris}
 &    E\Big[ z(Z_1, Z_2, Z_3, Z_4) - z(Z'_1, Z_2, Z_3, Z_4) \Big] \nonumber \\
      & = E\Big[\int_{\mathbb{R}^2} \Big(I(Z_1 > u) - I(Z'_1 > u) \Big)\Big(I(Z_2 > v) - I(Z'_2 > v) \Big)z''(u, v, Z_3, Z_4)dvdu\Big] + 0.
 \end{align}
As $h^{(2)}(. , .)$ is defined in  $(\ref{h2defn})$ is a degenerate kernel $(\int_{\mathbb{R}}h^{(2)}(x, y)dF(y) = 0, \text{ for all } x \in \mathbb{R})$,  $E(h^{(2)}(Z'_1, Z_2)h^{(2)}(Z_3,Z_4)) =0$ and $E(h^{(2)}(Z'_1, Z'_2) h^{(2)}(Z_3,Z_4)) =0$. Hence, $(\ref{eq1chap1chris})$ follows. Solving similarly as $(\ref{eq1chap1chris})$, we get
 \begin{align*} 
  &E\Big[ z(Z_1, Z_2, Z_3, Z_4) - z(Z'_1, Z_2, Z_3, Z_4) \Big] \nonumber \\
     & =  E\Big[\int_{\mathbb{R}^4} \Big(I(Z_1 > u) - I(Z'_1 > u) \Big)\Big(I(Z_2 > v) - I(Z'_2 > v) \Big)\Big(I(Z_3 > x) - I(Z'_3 > x) \Big)\nonumber \\ &\hspace{1in}\Big(I(Z_4 > y) - I(Z'_4 > y) \Big)z''''(u, v, x, y)dudvdxdy\Big].
 \end{align*}
\end{proof}

In the following, assume $\lbrace{X_i', i \ge 1}\rbrace$ to be a sequence of random variables independent of ${X_i, i \ge 1}$ such that ${X_i', i \ge 1}$ are  $i.i.d.$ with $F$ as the marginal distribution function of $X_1'$.

 \begin{lemma} \label{lemma2chap1}
Let the functions $h^{(2)}(. , .)$, $z(., ., ., .)$ and $z''''(., ., ., .)$ be as defined in Lemma $\ref{chris}$. Assume the following holds.
\begin{itemize}
\item[(i)]  $E\Big(h^{(2)}(X_1, X_{1+k})\Big)^{2}< \infty$, for all $k \in \mathbb{N}$;
\item[(ii)] For any $0 < C_1 < \infty$,  $\underset{u, v, x, y \in [-C_1, C_1]}{sup} ||z''''(u, v, x, y)||$ $ = CC_1^{b}$, for some $b \in \mathbb{N}\cup \{0\}$;
\item[(iii)] For the chosen $C_1$, define functions $f_1(.)$, and $f_2(.)$ as follows.
\begin{align*}
& f_1(x) = xI(|x| \le C_1) +C_1I(x > C_1) - C_1I(x < - C_1) \text{ and}\\
& f_2(x) =(x-C_1)I(x > C_1) + (x+C_1)I(x < -C_1).
\end{align*}
Then, for  all $k_1, k_3 \in \mathbb{N}$, $k_2 \in \mathbb{N}\cup \{0\}$, and some $\delta > 0$,
\begin{align*}
\sum_{\underset{(i, j,p,q) \neq (1, 1, 1, 1)}{{i,j,p,q=1,2}}}E\Big|\Big(h^{(2)}(f_{i}(Z_1), f_{j}(Z_{1+k_1}))\Big)\Big(h^{(2)}(f_p(Z_{1+k_2}), f_q(Z_{1+k_3}))\Big)\Big| \le \frac{C}{C_1^{\delta}}, 
\end{align*}
where $Z_m= X_m \text{ or } X'_m$, for $m = 1, 1+k_1,1+k_2, 1+k_3$;
\item[(iv)] $\sum_{j=1}^{\infty}Cov(X_1, X_j)^{\delta/(3(p+\delta))} < \infty$ for $p = b+2$, where b and $\delta$ are defined in $(ii)$ and $(iii)$, respectively.
\end{itemize}
 Then, as $n \to \infty$,
\begin{equation}\label{eq19chap1}
\sum_{1 \le i < j \le n} \sum_{1 \le k < l \le n} |E( h^{(2)}(X_i, X_j) h^{(2)}(X_k, X_l))| = o(n^3).
\end{equation}
\end{lemma}
\begin{proof}      Note that,  due to $(i)$,
\begin{equation}\label{eq8chap1}
\sum_{1 \le i < j \le n}  |E( h^{(2)}(X_i, X_j))^2| = o(n^3).
\end{equation}
Since $h^{(2)}(x, y)$ is a degenerate kernel, 
\begin{equation}
E( h^{(2)}(X_i', X_j) h^{(2)}(X_k, X_l)) = 0.
\end{equation}
Using Lemma $\ref{chris}$, for $i, j, k, l$  all distinct,
\begin{align*}
&|E( h^{(2)}(X_i, X_j) h^{(2)}(X_k, X_l))| =  |E(z(X_i, X_j, X_k, X_l)-  z(X'_i, X_j, X_k, X_l))| \\
 & =\Big| E\Big[ \int_{\mathbb{R}^4}\Big(I(X_i > u_i) - I(X'_i > u_i) \Big)\Big(I(X_j > u_j) - I(X'_j > u_j) \Big)\Big(I(X_k > u_k) - I(X'_k > u_k) \Big)\nonumber \\ &\hspace{1in}\Big(I(X_l > u_l) - I(X'_l > u_l) \Big)z''''(u_i, u_j, u_k, u_l)du_idu_jdu_kdu_l\Big]\Big| \nonumber \\
 &  =  \Big|E\Big[ \int_{\mathbb{R}^4}z''''(u_i, u_j, u_k, u_l) \prod_{t = i, j, k, l}\Big(I(X_t > u_t) - I(X'_t > u_t) \Big)\prod_{t = i, j, k, l}du_t\Big]\Big|.
\end{align*}
 Choose a  $C_1> 0$. Define $A = [-C_1, C_1]^4$.
\begin{align*}
& |E( h^{(2)}(X_i, X_j) h^{(2)}(X_k, X_l))|  =  |E(z(X_i, X_j, X_k, X_l)-  z(X'_i, X_j, X_k, X_l))| \nonumber \\
  & \le \Big|E\Big[  \int_{A} z''''(u_i, u_j, u_k, u_l)\prod_{t = i, j, k, l}\Big(I(X_t > u_t) - I(X'_t > u_t) \Big)\prod_{t = i, j, k, l}du_t\Big]\Big|  \nonumber \\
   & \hspace{0.5in}\ + \Big|E\Big[  \int_{A^c} z''''(u_i, u_j, u_k, u_l)\prod_{t = i, j, k, l}\Big(I(X_t > u_t) - I(X'_t > u_t) \Big)\prod_{t = i, j, k, l}du_t\Big]\Big| \nonumber \\
   & = |I_1|+ |I_2|.
\end{align*}
 \begin{align}
       & \Big|I_1\Big| =   \Big|\int_{A}z''''(u_i, u_j, u_k, u_l) E\Big[\prod_{t = i, j, k, l}\Big(I(X_t > u_t) - I(X'_t > u_t) \Big)\Big] \prod_{t = i, j, k, l}du_t\Big| \nonumber \\
       & \le C \underset{u, v, x, y \in [-C_1, C_1]}{sup} |z''''(u, v, x, y)|\int_{A}\Big|E\Big[\prod_{t = i, j, k, l}\Big(I(X_t > u_t) - I(X'_t > u_t) \Big)\Big]\Big|\prod_{t = i, j, k, l}du_t \nonumber \\
& \le CC_1^b C_1^2\Big|\int_{[-C_1, C_1]^2} \Big[\sum_{t = j, k, l}P(X_i > u_i, X_t > u_t) - P(X_i > u_i)P( X_t > u_t)du_idu_t \Big]\Big| \label{lebowitz123}\\
& \le CC_1^p\Big|\int_{\mathbb{R}^2} \Big[\sum_{t = j, k, l}P(X_i > u_i, X_t > u_t) - P(X_i > u_i)P( X_t > u_t)du_idu_t\Big] \Big| \nonumber \\
& = CC_1^p [Cov(X_i, X_j) + Cov(X_i, X_k) + Cov(X_i, X_l) ] \label{lebowitz1234}.
 \end{align}
The inequality in $(\ref{lebowitz123})$ follows from Lebowitz's inequality, Lemma $\ref{lebowitz}$. The equality in $(\ref{lebowitz1234})$ follows using the Newman's inequality, Lemma $\ref{covnew1}$.  Hence,
\begin{equation}\label{eq4chap1}
 |I_1| \le CC_1^p[Cov(X_i, X_j) + Cov(X_i, X_k) + Cov(X_i, X_l)]= CC_1^p[\sum_{t_1 = j, k, l}Cov(X_i, X_{t_1})].
 \end{equation}
 Similarly, it can be shown
\begin{align}
&   |I_1|   \le CC_1^p(Cov(X_j, X_i) + Cov(X_j, X_k)  + Cov(X_j, X_l) ) = CC_1^p[\sum_{t_2 = i, k, l}Cov(X_j, X_{t_2})]. \label{eq2chap1} \\
& |I_1|  \le CC_1^p(Cov(X_k, X_i) + Cov(X_k, X_j)  + Cov(X_k, X_l) )= CC_1^p[\sum_{t_3 = i, j, l}Cov(X_k, X_{t_3})].\label{eq3chap1}
\end{align}
Combining $(\ref{eq4chap1})$ - $(\ref{eq3chap1})$, we get the following.
\begin{align}\label{eq11chap1}
 |I_1|  \le  CC_1^pT^{1/3},
\end{align}
where $T = [\sum_{t_1 = j, k, l}Cov(X_i, X_{t_1})]\times[\sum_{t_2 = i, k, l}Cov(X_j, X_{t_2})]\times[\sum_{t_3 = i, j, l}Cov(X_k, X_{t_3})]$.  Next, 
\begin{align*}
 I_2 =   &E\Bigg[  \int_{A^c}z''''(u_i, u_j, u_k, u_l) \prod_{t = i, j, k, l}\Big(I(X_t > u_t) - I(X'_t > u_t) \Big)\prod_{t = i, j, k, l}du_t\Bigg].
\end{align*}
There are several combinations possible under $A^c$. We solve the integral for a  possible combination under $A^c$, say $B = \{ |u_i| > C_1, |u_j| \le C_1, |u_k| \le C_1, |u_l| \le C_1\}$.  Calculations for other combinations will follow similarly.

For a differentiable function $g(.)$ with derivative $g'(.)$, and random variables $X$ and $X'$, we have
\begin{align}
& \int_{|v| \le C_1} \Big(I(X > v) - I(X' > v) \Big)g'(v)dv = g(f_1(X)) - g(f_1(X')) \nonumber \\
& \int_{|v| > C_1} \Big(I(X > v) - I(X' > v) \Big)g'(v)dv = g(f_2(X)) - g(f_2(X')).
\end{align}
Using this, we have 
\begin{align*}
&E\Bigg[  \int_{B}z''''(u_i, u_j, u_k, u_l) \prod_{t = i, j, k, l}\Big(I(X_t > u_t) - I(X'_t > u_t) \Big)\prod_{t = i, j, k, l}du_t\Bigg] \nonumber \\
& = E\Bigg[ \Big( h^{(2)}(f_2(X_i), f_1(X_j)) -  h^{(2)}(f_2(X_i), f_1(X_j')) - h^{(2)}(f_2(X'_i), f_1(X_j)) + h^{(2)}(f_2(X'_i), f_1(X_j'))\Big) \times \nonumber \\ 
& \hspace{0.4in} \Big( h^{(2)}(f_1(X_k), f_1(X_l)) -  h^{(2)}(f_1(X_k), f_1(X_l')) - h^{(2)}(f_1(X'_k), f_1(X_l)) + h^{(2)}(f_1(X'_k), f_1(X_l'))\Big)\Bigg].
\end{align*}
From $(iii)$, we get
\begin{align*}
&\Bigg|E\Bigg[  \int_{B}z''''(u_i, u_j, u_k, u_l) \prod_{t = i, j, k, l}\Big(I(X_t > u_t) - I(X'_t > u_t) \Big)\prod_{t = i, j, k, l}du_t\Bigg] \Bigg|\le \frac{C}{C_1^{\delta}}.\label{eq14chap1}
 \end{align*}
 Similarly, solving for other combinations under $A^{c}$, we get 
 \begin{equation}\label{eq15chap1}
 |I_2| \le \frac{C}{C_1^{\delta}}.
 \end{equation}
From $(\ref{eq11chap1})$ and $(\ref{eq15chap1})$, for $i, j, k, l$  all distinct, we have
\begin{equation}\label{eq1578chap1}
 |E( h^{(2)}(X_i, X_j) h^{(2)}(X_k, X_l))| \le CC_1^pT^{1/3} + \frac{C}{C_1^{\delta}}.
\end{equation}
Choosing $C_1 = T^{-1/(3(p+\delta))}$, we get,
\begin{equation}\label{eq17chap1}
 |E( h^{(2)}(X_i, X_j) h^{(2)}(X_k, X_l))| \le CT^{\delta/(3(p+\delta))}.
\end{equation}
Next, consider a case when there are 3 distinct indices in $i, j, k, l$. Let $j = k$, then solving as in  $(\ref{eq4chap1})$ and $(\ref{eq15chap1})$
\begin{equation}\label{eqblahchap1}
  |E( h^{(2)}(X_i, X_j) h^{(2)}(X_j, X_l))|   \le CC_1^p[Cov(X_i, X_j)  + Cov(X_i, X_l)] + \frac{C}{C_1^{\delta}}
\end{equation}
Choosing $C_1 = T^{-1/(p+\delta)}$, we get,
\begin{equation}\label{eq18chap1}
 |E( h^{(2)}(X_i, X_j) h^{(2)}(X_j, X_l))| \le C[Cov(X_i, X_j)^{\delta/(\delta+p)}  + Cov(X_i, X_l)^{\delta/(\delta+p)} ]  
\end{equation}
Using $(\ref{eq8chap1})$, $(\ref{eq17chap1})$, $(\ref{eq18chap1})$, stationarity of $X_j's$ and $(iv)$, $(\ref{eq19chap1})$ follows.
 \end{proof}
 
 \begin{remark} The condition $(iii)$ of Lemma $\ref{lemma2chap1}$ may seem cumbersome, but  in general can be shown to be true under restrictions on the moments of $X_1$, as seen in Theorem $\ref{mu2chap1}$.
  \end{remark}
If the random variables  ${X_n, n \ge 1}$  are uniformly bounded (as is often seen in  applications in the reliability and survival analysis) then we can use the following lemma. The assumptions on the covariance structure become less restrictive. 
\begin{lemma} \label{lemma1chap1}
Let the functions $h^{(2)}(., .)$, $z(., ., ., .)$ and $z''''(., ., ., .)$ be as defined in Lemma $\ref{chris}$. Assume the following holds.
\begin{itemize}
\item[(a)] $X_n, n \ge 1$ are uniformly bounded, $i.e.$ $P(|X_1| < C_1) = 1$, for some $C_1 > 0$; \item[(b)] $\sum_{j=1}^\infty Cov(X_1, X_j)^{1/3}$ $<$ $\infty$;
\item[(c)] $|h^{(2)}(x, y)|$ is bounded for all $x, y \in [-C_1, C_1]$, where $h^{(2)}(x, y)$ is as defined in $(\ref{h2defn})$;
\item[(d)]  $ \underset{u, v, x, y \in [-C_1, C_1]}{sup} ||z''''(u, v, x, y)||$ is  bounded.
\end{itemize}
 Then, as $n \to \infty$,
\begin{equation}\label{eq7chap1}
\sum_{1 \le i < j \le n} \sum_{1 \le k < l \le n} |E( h^{(2)}(X_i, X_j) h^{(2)}(X_k, X_l))| = o(n^3).
\end{equation}
\end{lemma}
\begin{proof} The proof follows similarly as Lemma $\ref{lemma2chap1}$. The second term $\frac{C}{C_1^{\delta}}$, for some $\delta> 0$ in the upper bounds obtained in  inequalities $(\ref{eq1578chap1})$ and $(\ref{eqblahchap1})$ would not be needed.
\end{proof}

\begin{lemma}\label{lemma3chap1}
Let $U_n$ be the U-statistic based on a symmetric kernel $\rho(x, y)$, such that the corresponding kernel $h^{(2)}(., .)$ (defined in (\ref{h2defn})) satisfies the conditions of Lemma $\ref{lemma2chap1}$ or Lemma $\ref{lemma1chap1}$. Let  $0 < \sigma_1^2$ $=$ $Var(\rho_1(X_1))$  $<$ $\infty$. 
Define $\sigma_{1j}= Cov( \rho_1(X_1), \rho_1(X_{1+j}))$. Assume $\sum_{j=1}^{\infty}|\sigma_{1j}| < \infty$.
Then,
\begin{equation}\label{varchap1eq20}
Var(U_n) = \frac{4\sigma^2_U}{n} + o\Big(\frac{1}{n}\Big), \text{ where } \sigma^2_U = \sigma^2_1 + 2\sum_{j=1}^{\infty}\sigma_{1j}.
\end{equation}
\end{lemma}
\begin{proof} The proof follows similarly as Lemma $3.2$ of \cite{Garg2015209}. By Hoeffding's decomposition, $Var(U_n) = 4Var(H^{(1)}_n) + Var(H^{(2)}_n) + 4Cov(H^{(1)}_n, H^{(2)}_n)$. As  $\sum_{j=1}^{\infty}|\sigma_{1j}| < \infty$,
\begin{equation}
Var(H^{(1)}_n) = \frac{1}{n}(\sigma^{2}_{1} +2\sum_{j=1}^{\infty}\sigma_{1j})
                          + o\Big(\frac{1}{n}\Big).
\end{equation}
 Also $Var(H^{(2)}_{n}) \le E(H^{(2)}_{n})^2  = o({1}/{n})$ (from Lemma $\ref{lemma2chap1}$ or Lemma $\ref{lemma1chap1}$).
Using Cauchy-Schwarz inequality we have,  $|Cov(H^{(1)}_n, H^{(2)}_n)| \le  o({1}/{n})$. Hence, we get $(\ref{varchap1eq20})$.
\end{proof}

The following gives the central limit theorem for a non-degenerate U-statistic based on a stationary sequence of associated observations with  a  kernel, $\rho$, of degree 2.

\begin{theorem}\label{theorem1chap1}
Suppose the conditions of Lemma $\ref{lemma3chap1}$ hold. Further, assume $\sigma^2_U > 0$, where $\sigma^2_U$ is defined by $(\ref{varchap1eq20})$. If there exists a function $\tilde{\rho}_1(\cdot)$ such that ${\rho}_1$ $\ll$ $\tilde{\rho}_1$ and
\begin{equation}\label{eq21chap1}
\sum_{j=1}^{\infty}Cov(\tilde{\rho}_1(X_1), \tilde{\rho}_1(X_{j})) < \infty,
\end{equation}
then
\begin{equation}
\frac{\sqrt{n}( U_n - \theta)}{2\sigma_U} \xrightarrow {\mathcal{ L}} N(0, 1) \:\: \text{as}\: \:{n \to \infty}.
\end{equation}
\end{theorem}
\begin{proof} Using  Hoeffding's decomposition for $U_n$,
\begin{equation}\label{hoeffd}
\frac{\sqrt{n}( U_n - \theta)}{2\sigma_U} = n^{{-1}/{2}}\sum_{j=1}^{n}\frac{h^{(1)}(X_j)}{\sigma_U} + \frac{\sqrt{n}H^{(2)}_n}{\sigma_U}.
\end{equation}
In addition, $nE({H^{(2)}_n})^2 \xrightarrow \: 0 \: \text{as}\: \:{n \to \infty}$, from Lemma $\ref{lemma2chap1}$ or Lemma $\ref{lemma1chap1}$. Hence,
\begin{equation}\label{On8}
\frac{\sqrt{n}H^{(2)}_n}{\sigma_U}\xrightarrow {p}\: 0 \: \text{as}\: \:{n \to \infty}.
\end{equation}
From  Lemma \ref{cltfass}, we get that,
\begin{equation}\label{normU}
n^{{-1}/{2}}\sum_{j=1}^{n}\frac{h^{(1)}(X_j)}{\sigma_U} \xrightarrow {\mathcal{ L}} N(0, 1) \:\: \text{as}\: \:{n \to \infty}.
\end{equation}
Relations {$(\ref{hoeffd})$},  {$(\ref{On8})$} and  $(\ref{normU})$ prove the theorem.
\end{proof}

\begin{remark} \label{remarkUchap1general}The above results can be easily extended to a U-statistic $U_n$ based on a symmetric kernel $\rho(x_1,\dots, x_k)$ which is of a finite degree $k > 2$.  Let $\lbrace{X_i', i \ge 1}\rbrace$ be a sequence of random variables as defined earlier. Assume the following.
\begin{itemize}
\item[(i)] $E(\rho(X_{i_1}, \cdots, X_{i_k}))^{2} < \infty$, for all $i_j \in \mathbb{N}$, $j = 1, \cdots, k$, such that $i_1 < \cdots < i_k$. 
\item[(ii)] For any $0 < C_1 < \infty$, $z(x_1, \cdots, x_k, x'_1, \cdots, x'_k) = \rho({x_1},  \dots, {x_k})\rho({x'_1}, \dots, {x'_k})$ has  a bounded derivative $\frac{\partial^{2k}z}{\prod_{i = 1}^k\partial x_i \partial x_i'}$  for all $x_i, x'_i \in [-C_1, C_1]$, $i = 1, 2, \cdots, k$, and the bound is $CC_1^b$, for some $b \in \mathbb{N}\cup \{0\}$;
\item[(iii)] For some  $\delta > 0$ and all $j_p \in \mathbb{N}$, $p = 1, \cdots, 2k$, such that $ j_1 < \cdots < j_k$ and  $ j_{k+1} < \cdots < j_{2k}$, 
\begin{align*}
& \sum_{\underset{(i_1, \cdots,i_{2k}) \neq (1, \cdots, 1)}{{i_1,\cdots,i_{2k}=1,2}}}E\Big|z(f_{i_1}(Z_{j_1}), f_{i_2}(Z_{j_2}), \cdots, f_{i_k}(Z_{j_k}), f_{i_{k+1}}(Z_{j_{k+1}}),  \cdots, f_{i_{2k}}(Z_{j_{2k}}))\Big| \le \frac{C}{C_1^{\delta}}, 
\end{align*}
where $Z_m = X_m \text{ or } X'_m$, $m = j_1, j_2, \cdots, j_{2k}$;
 \item[(iv)] $\sum_{j=1}^{\infty}Cov(X_1, X_j)^{\delta/(3(p+\delta))} < \infty$ for $p = b+2k-2$, where $b$ and $\delta$ are respectively defined in $(ii)$ and $(iii)$ respectively;
 \item[(v)] $0 < \sigma^2_1 < \infty$ and $\sum_{j=1}^{\infty} |\sigma_{1j}| < \infty$, with $\sigma_U^2 > 0$; and
 \item[(vi)] There exists a function $\tilde{\rho}_1$,  such that ${\rho}_1$ $\ll$ $\tilde{\rho}_1$ and $(\ref{eq21chap1})$ holds. \end{itemize}
Then,
\begin{equation}
\frac{\sqrt{n}( U_n - \theta)}{k\sigma_U} \xrightarrow {\mathcal{ L}} N(0, 1) \:\: \text{as}\: \:{n \to \infty}. 
\end{equation}
\end{remark}

\begin{remark} Using Lemma $\ref{covnew3}$, we get $\sigma_1^2 \le Var(\tilde{\rho}_1(X_1))$and $|\sigma_{1j}|  \le  C \: Cov(\tilde{\rho}_1(X_1), \tilde{\rho}_1(X_{j}))$. If $\rho_1$ is monotonic, then $\{\rho_1(X_n), n \ge 1\}$ is a sequence of stationary associated random variables and one can take $\tilde{\rho}_1 \equiv \rho_1$.
\end{remark}
\begin{remark} If $\tilde{\rho}_1$ has a bounded derivative, then using  Lemma $\ref{covnew1}$, if  $\sum_{j=1}^{\infty}Cov(X_1, X_j) \: < \: \infty$, then  $(\ref{eq21chap1})$ holds.
\end{remark}
\begin{remark} Let $\tilde{\rho}_1(x)$ $=$ $cx$ for some constant $c > 0$ for all $x \in \mathbb{R}$. If ${\rho}_1$ $\ll$ $\tilde{\rho}_1$ then ${\rho}_1$ is  a Lipschitz function. A sufficient condition for ${\rho}_1$ to be a  Lipschitz function is that it should have a bounded derivative.
\end{remark}
\begin{remark}
If ${\rho}_1(x)$ is a function of bounded variation, then there exist two increasing functions $U_1(x)$ and $U_2(x)$ such that ${\rho}_1(x)= U_1(x) -  U_2(x)$ for all $x \in \mathbb{R}$. ${\rho}_1 \ll \tilde{\rho}_1$ by taking $\tilde{\rho}_1(x)= U_1(x) + U_2(x)$.  
\end{remark}

\section{Applications}\label{section4}
\setlength\parindent{24pt}
We use the results derived in section $\ref{section3}$ to get the limiting distribution of estimators of second, third and fourth central  moments, when the underlying sample is from ${X_n, n \ge 1}$. As given earlier,  $\{X_n, n \ge 1 \}$ is a sequence of stationary associated random variables, with $F$ as the marginal distribution function of $X_1$. We  also discuss estimators of skewness and kurtosis based on the estimators of moments. 

Define, $\mu = E(X_1)$, $\mu'_k = E(X_1^k)$ and $\mu_k = E((X_1-\mu)^k)$, $k \ge 2$. Assume that $\mu_k$ and $\mu'_k$ exist  for all $k = 2, \cdots, 8$, and $C$ is a generic positive constant in the sequel.

\subsection{Estimator of the variance/ second central moment}
A measure of variability  is the variance $\mu_2 = E((X_1-\mu)^2)$.

Given the sample ${X_j, 1 \le j \le n }$ from $F$, a U-statistic, $\hat{\mu}_{2n}$, estimating $\mu_2$ is
\begin{align}
\hat{\mu}_{2n}  = \frac{1}{(n-1)}\sum_{1 \le i \le n} (X_i - \bar{X}_n)^2  = \frac{2}{n(n-1)}\sum_{1 \le i <j \le n} \rho(X_i, X_j),
\end{align}
where $\bar{X}_n$ $=$ $\frac{\sum_{i=1}^n X_i}{n}$ and the kernel $\rho(x, y) = \frac{(x - y)^2}{2}$, with $\rho_1(x) = ({(x- \mu)^2 + \mu_2})/{2}$.

We next discuss the limiting distribution of $\hat{\mu}_{2n}$,  when the sample is from $\{X_n, n \ge 1\}$.

\begin{theorem} \label{mu2chap1}Let  $E|X_1^{\gamma}| < \infty$, where $\gamma = {(4+\delta)}$ for some $\delta >0$. Further, if
\begin{equation}\label{eq25chap1}
 \sum_{j=1}^{\infty}Cov(X_1, X_j)^{\frac{\delta}{{3(2+\delta)}}} \: < \: \infty,
\end{equation}
then
\begin{equation}\label{eq22chap1}
\frac{\sqrt{n}(\hat{\mu}_{2n} - \mu_2)}{2\sigma_{\hat{\mu}_{2n}}}  \xrightarrow {\mathcal{ L}} N(0, 1) \:\: \text{as}\:\: {n \to \infty},
\end{equation}
where ${\sigma_{\hat{\mu}_{2n}}} = Var(\rho_1(X_1)) + 2\sum_{j=2}^{\infty} Cov(\rho_1(X_1), \rho_1(X_j))$ and $Var(\rho_1(X_1))= \mu_4 - ({\mu_2})^2$.
 
If ${X_n, n \ge 1}$ are uniformly bounded, then $(\ref{eq22chap1})$ holds under
  \begin{equation}\label{eq23chap1}
 \sum_{j=1}^{\infty}Cov(X_1, X_j)^{1/3} \: < \: \infty.
 \end{equation}
\begin{proof}  Assume ${X_n, n \ge 1}$ are uniformly bounded, i.e. $P(|X_1| < C_1) = 1$, for some $C_1 > 0$.  As $z(x, y, x', y') = h^{(2)}(x, y)h^{(2)}(x', y') =(x-\mu)(y-\mu)(x'-\mu)(y'-\mu)$ and $z''''(x, y, x', y') = 1$, the conditions of  Lemma $\ref{lemma1chap1}$ are true under $(\ref{eq23chap1})$.  As $\rho_1(x)$ $\ll$ $Cx$, for all $x \in  [-C_1, C_1]$, conditions of  Theorem $\ref{theorem1chap1}$ are satisfied, and we get $(\ref{eq22chap1})$.

The above can be easily extended to random variables which are not uniformly bounded, under the conditions of Lemma $\ref{lemma2chap1}$. Without loss of generality assume $\mu = 0$. The conditions $(i) - (iv)$ of Lemma $\ref{lemma2chap1}$ are satisfied under the given assumptions, as discussed in the following.
\begin{itemize}
\item[(i)]  Assuming $E|X_1|^{4 + \delta}$ $<$ $\infty$,  $(i)$ of Lemma $\ref{lemma2chap1}$ is  true. As $h^{(2)}(x, y) =xy$,
\begin{align*}
E(h^{(2)}(X_1, X_{1+k}))^2 = E(X_1 X_{1+k})^{2} \le CE(X_1)^{4} \le C(E|X_1|^{4 + \delta})^{4/(4+\delta)}< \infty, \text{ for all $k \in \mathbb{N}$}.
 \end{align*}
 \item[(ii)] $b= 0$, as  $z''''(x, y, x', y') = 1$, for all $x, y, x', y'$ $\in$ $\mathbb{R}$.
 \item[(iii)] Assuming $E|X_1|^{4 + \delta}$ $<$ $\infty$,  $(iii)$ of Lemma $\ref{lemma2chap1}$ is  true. Consider
  \begin{align*}
& \Big|E\Big(f_2(X_i)f_2(X_j)f_1(X_k)f_1(X_l)\Big)  \Big| \le C_1^2 \Big|E\Big(|X_i| |X_j| I(|X_i| > C_1) I(|X_j| > C_1)\Big)\Big| \\
& \le  C_1^2 E\Big(|X_i| |X_j| \frac{|X_i|^{1 + \delta/2} }{ C_1^{1 + \delta/2}}\frac{|X_j|^{1 + \delta/2}}{C_1^{1 + \delta/2}}\Big)\Big|  \le \frac{C}{C_1^\delta}.
    \end{align*}
Similarly, the bounds on the other terms can be obtained.
  \item[(iv)]  From $(ii)$, $p=2$ as $b=0$. The restriction put on the covariance structure is therefore $\sum_{j=1}^{\infty}Cov(X_1, X_j)^{\frac{\delta}{{3(2+\delta)}}} \: < \: \infty$. 
   \end{itemize} 

Define
\begin{equation*}
\tilde{\rho}_1(x) = \frac{x^2I(x \ge 0)   - x^2I(x < 0)}{2},
\end{equation*}
where $I(.)$ denotes the indicator function. Then $\rho_1(x)$ $\ll$ $\tilde{\rho}_1(x)$. It can be shown that for any $0 < C_3 <\infty$,
\begin{equation} \label{eq26chap1}
  Cov(\tilde{\rho}_1(X_i), \tilde{\rho}_1(X_j)) \le \frac{C}{C_3^\delta} + C C_3^2 Cov(X_i, X_j).
  \end{equation}
  Choosing, $C_3 = Cov(X_i, X_j)^{-1/{(2+\delta)}}$ in $(\ref{eq26chap1})$,
  \begin{equation*}
   Cov(\tilde{\rho}_1(X_i), \tilde{\rho}_1(X_j)) \le CCov(X_i, X_j)^{\delta/{(2+\delta)}}.
 \end{equation*}
Using $\sum_{j=2}^{\infty} Cov(X_1, X_j)^{\frac{\delta}{2+\delta}} < \infty$ (which follows from $(\ref{eq25chap1})$) and Theorem $\ref{theorem1chap1}$, we get $(\ref{eq22chap1})$.
\end{proof}
\end{theorem}

 \begin{remark}   The limiting distribution of the variance estimator, $\hat{\mu}_{2, n}$,  when the underlying sample is from ${X_n, n \ge 1}$, can also be obtained  under different set of assumptions.  
  \begin{itemize}
  \item[(i)]To apply the techniques of  \cite{beutner2012, beutner2014} to obtain a central limit theorem for  $\hat{\mu}_{2, n}$, a result on weak convergence of weighted empirical distribution process of associated sequences is needed. This result  can be found in \cite{shao1996}. Our results are  valid for non-degenerate U-statistics only, while the technique of \cite{beutner2014}  can also be used when the U-statistics considered
are degenerate.\\ The following is given in Example $3.8$ of \cite{beutner2012}.  \\Let $\{X_i, i \ge 1\}$ be a stationary, associated sequence with the marginal distribution function $F$ and $Cov(X_1, X_n) = O(n^{-\nu - \epsilon})$ for some $\nu \ge (3 +\sqrt{33})/2 \approx 4.372281$ and $\epsilon > 0$. Then, whenever $F$ has a finite  $\gamma$-moment for some  $\gamma> \frac{4\nu}{\nu - 3}$, a central limit theorem for $\hat{\mu}_{2, n}$ holds.
   \item[(ii)] The kernel $\rho(x_1, x_2) = (x_1 - x_2)^2/2$ is a function of bounded Hardy-Krause variation on any bounded rectangle $[-C_1, C_1]^2$  where $0 < C_1 < \infty$. This can be shown as following.  We can write $(x_1 - x_2)^2/2$ as a difference of  two component-wise nondecreasing bounded functions $g_1(x_1, x_2) = ((x_1 + C_1)^2 + (x_2 + C_1)^2)/2$ and $g_2(x_1, x_2) = (x_1 + C_1)(x_2 + C_1)$, where for $x_1 \le y_1$ and $x_2 \le y_2$, $g_i(x_1, x_2) - g_i(x_1, y_2) - g_i(y_1, x_2) + g_i(y_1, y_2) \ge 0$, for $i = 1, 2$. Using the results of \cite{Garg2015209}, a central limit theorem for the variance estimator holds if $\sum_{j=1}^{\infty}Cov(X_1, X_j)^{\nu} < \infty$, for some $0 < \nu < 1/6$, when $\{X_n, n \in \mathbb{N}\}$ is uniformly bounded. The covariance structure assumed in Theorem $\ref{mu2chap1}$ is less restrictive.  Also, the condition of differentiability is relatively easier to verify. 
 \end{itemize}
\end{remark}

 \begin{remark}
 The results of \cite{Dewan20029} and its corrigendum \cite{Dewan2015147} cannot be used to obtain the central limit theorem for the variance estimator. They only consider the  case when the differentiable kernels are component-wise monotonic. Moreover, their results are based on a very restrictive set of conditions. They assume that
there exists a non-negative function $r(k)$ satisfying $\sum_{k =0}^{\infty} r(k) < \infty$, such that for all $i, j, k, l$,
\begin{equation*}
|Cov(h^{(2)}(X_i, X_j), h^{(2)}(X_k, X_l))| \le r(max(|i-k|, |j - l|)), 
\end{equation*}
They do not discuss  examples of functions $r(k)$ that satisfy the above inequality. The cases when kernels have a degree greater than 2 are also not considered.  
\end{remark}

\subsection{Estimator for the third central moment }
 Let the sample $X_j,  1 \le j \le n$ be from $F$. A U-statistic, $\hat{\mu}_{3n}$, based on  $ X_j,  1 \le j \le n$ estimating the third central moment $\mu_3$ is
\begin{align*}
\hat{\mu}_{3n}  = \frac{1}{{n \choose3}} \sum_{1 \le i < j < k \le n} u( X_i, X_j, X_k), \text{ where}
\end{align*}
\begin{align*}
u( x_1, x_2, x_3)  = & \sum_{i=1}^3\frac{x_i^3}{3}  - \sum_{(i, j, k) \in \pi(1, 2, 3)}\frac{x_i^2( x_j+ x_k)}{2} +  2[x_1x_2x_3],
\end{align*}
and  $\pi(1, 2, 3) = ((1, 2, 3), (2, 1, 3), (3, 1, 2))$. Here, $u_1(x) =   ({2\mu_3 + (x-\mu)^3})/{3} - \mu_2(x-\mu)$.

  \begin{theorem}\label{mu3chap1} Let  $E|X_1^{\gamma}| < \infty$, where $\gamma = {(6+\delta)}$ for some $\delta >0$. Further, if  \\$\sum_{j=1}^{\infty}Cov(X_1, X_j)^{\frac{\delta}{{3(4+\delta)}}} \: < \: \infty$, then
\begin{equation}\label{eq27chap1}
\frac{\sqrt{n}(\hat{\mu}_{3n} - \mu_3)}{3\sigma_{\hat{\mu}_3}}  \xrightarrow {\mathcal{ L}} N(0, 1) \:\: \text{as}\:\: {n \to \infty},
\end{equation}
where ${\sigma_{\hat{\mu}_3}^2} = Var(u_1(X_1)) + 2\sum_{j=2}^{\infty} Cov(u_1(X_1), u_1(X_j))$ and  $Var(u_1(X_1))= \mu_6- \mu_3^2$.

If ${X_n, n \ge 1}$ are uniformly bounded, then $(\ref{eq27chap1})$ is true under $\sum_{j=1}^{\infty}Cov(X_1, X_j)^{1/3} \: < \: \infty$.

 \end{theorem}
 \begin{proof} Proof follows similarly as the proof of Theorem $\ref{mu2chap1}$. The conditions are obtained by putting $k =3$, $b= 0$ and $p= 4$ in Remark $\ref{remarkUchap1general}$.
  \end{proof}
 \subsection{Estimator for the fourth central moment }
 Let the sample $X_j,  1 \le j \le n$ be from $F$. A U-statistic, $\hat{\mu}_{4n}$,  based on this sample estimating the fourth central moment $\mu_4$ is
\begin{align*}
\hat{\mu}_{4n}= \frac{1}{{n \choose4}} \sum_{1 \le i < j < k < l \le n} v( X_i, X_j, X_k, X_l), \text{  where}
\end{align*}
\begin{align*}
 v( x_1, x_2, x_3, x_4)  = & \sum_{i=1}^4\frac{x_i^4}{4}  + \sum_{(i, j, k, l) \in \pi(1, 2, 3, 4)}\frac{x_i^2( x_jx_k + x_kx_l +x_jx_l)}{2}  \\ 
 &  - \sum_{(i, j, k, l) \in \pi(1, 2, 3, 4)} \frac{x_i^3(x_j+x_k+ x_l)}{3}   - 3[x_1x_2x_3x_4],
\end{align*}
and $\pi(1, 2, 3, 4) = ((1, 2, 3, 4), (2, 1, 3, 4), (3, 1, 2, 4), (4, 1, 2, 3))$. 

Here, $v_1(x) =   ({3\mu_4 + (x-\mu)^4})/{4} - \mu_3(x-\mu)$.

 \begin{theorem} \label{mu4chap1}Let  $E|X_1^{\gamma}| < \infty$, where $\gamma = {(8+\delta)}$ for some $\delta >0$. Further, if  \\$\sum_{j=1}^{\infty}Cov(X_1, X_j)^{\frac{\delta}{{3(6+\delta)}}} \: < \: \infty$, then
\begin{equation}\label{eq28chap1}
\frac{\sqrt{n}(\hat{\mu}_{4n} - \mu_4)}{4\sigma_{\hat{\mu}_4}}  \xrightarrow {\mathcal{ L}} N(0, 1) \:\: \text{as}\:\: {n \to \infty},
\end{equation}
where ${\sigma_{\hat{\mu}_4}^2} = Var(v_1(X_1)) + 2\sum_{j=2}^{\infty} Cov(v_1(X_1), v_1(X_j))$ and  $Var(v_1(X_1))= \mu_8 - (\mu_4)^2$.
 
If ${X_n, n \ge 1}$ are uniformly bounded, then $(\ref{eq28chap1})$ is true under $\sum_{j=1}^{\infty}Cov(X_1, X_j)^{1/3} \: < \: \infty$.

 \end{theorem}
  \begin{proof} Proof follows similarly as the proof of Theorem $\ref{mu2chap1}$. The conditions are obtained by putting $k =4$, $b= 0$  and $p= 6$ in Remark $\ref{remarkUchap1general}$.
  \end{proof}
  
 \begin{remark} The limiting distributions of U-statistic estimators of higher order central moments can be obtained similarly. In general, the $r^{th}$ central moment, $r \in \mathbb{N}/\{1\}$, is estimated by $\frac{\sum_{j=1}^n (X_j - \bar{X}_n)^r}{n}$. This and  the U-statistic estimator for the $r^{th}$ central moment will have the same limiting distribution.
 \end{remark}
 \subsection{ Estimation of skewness and kurtosis.}
 The coefficient of skewness and kurtosis  for the distribution function $F$ are defined as \\$ \tau = \frac{\mu_3}{(\mu_2)^{3 /2}} \text{ and } \kappa = \frac{\mu_4}{(\mu_2)^{2}}$,
  respectively. We can use 
  \begin{equation}\label{estichap2}
  \hat{\tau}_n= \frac{\hat{\mu}_{3n}}{(\hat{\mu}_{2n})^{3/2}} \text{ and }\hat{\kappa}_n= \frac{\hat{\mu}_{4n}}{(\hat{\mu}_{2n})^{2}}
  \end{equation}
  as  estimators of $\tau$ and $\kappa$, respectively.   
  
   Assume that  $\{X_n, n \in \mathbb{N} \}$ is a sequence of stationary associated random variables, with $F$ as the  marginal distribution function of $X_1$. Then, under the conditions of Theorems  $\ref{mu2chap1}$, $\ref{mu3chap1}$ and  $\ref{mu4chap1}$, it can be shown that $\hat{\tau}_n \xrightarrow{p} \tau \text{ and } \hat{\kappa}_n \xrightarrow{p} \kappa \text{ as } {n \to \infty}.$

\nocite{*}
\bibliographystyle{yearfirst}
\setlength{\bibsep}{0.0pt}
{\footnotesize\bibliography{references1}}

\begin{thebibliography}{16}
\providecommand{\natexlab}[1]{#1}
\expandafter\ifx\csname urlstyle\endcsname\relax
  \providecommand{\doi}[1]{doi:\discretionary{}{}{}#1}\else
  \providecommand{\doi}{doi:\discretionary{}{}{}\begingroup
  \urlstyle{rm}\Url}\fi

\bibitem[{Beutner and Zähle(2012)}]{beutner2012}
Beutner E. and Zähle H. (2012).
\newblock Deriving the asymptotic distribution of ${U}$- and ${V}$-statistics
  of dependent data using weighted empirical processes.
\newblock {\em Bernoulli\/}, 18(3):803--822.

\bibitem[{Beutner and Zähle(2014)}]{beutner2014}
Beutner E. and Zähle H. (2014).
\newblock Continuous mapping approach to the asymptotics of ${U}$- and
  ${V}$-statistics.
\newblock {\em Bernoulli\/}, 20(2):846--877.

\bibitem[{Bulinski and Shashkin(2007)}]{Bulinski:1701595}
Bulinski A. and Shashkin A. (2007).
\newblock {\em {Limit theorems for associated random fields and related
  systems}\/}.
\newblock Advanced Series on Statistical Science and Applied Probability. World
  Scientific, Singapore.

\bibitem[{Christofides and Vaggelatou(2004)}]{CHRISTOFIDES2004138}
Christofides T.C. and Vaggelatou E. (2004).
\newblock A connection between supermodular ordering and positive/negative
  association.
\newblock {\em J. Multivariate Anal.\/}, 88(1):138 -- 151.

\bibitem[{Dewan and Prakasa~Rao(2001)}]{dewanprakasarao2001}
Dewan I. and Prakasa~Rao B.L.S. (2001).
\newblock Asymptotic normality for {$U$}-statistics of associated random
  variables.
\newblock {\em J. Statist. Plann. Inference\/}, 97(2):201--225.

\bibitem[{Dewan and Prakasa~Rao(2002)}]{Dewan20029}
Dewan I. and Prakasa~Rao B.L.S. (2002).
\newblock Central limit theorem for {$U$}-statistics of associated random
  variables.
\newblock {\em Statist. Probab. Lett.\/}, 57(1):9 -- 15.

\bibitem[{Dewan and Prakasa~Rao(2015)}]{Dewan2015147}
Dewan I. and Prakasa~Rao B.L.S. (2015).
\newblock Corrigendum to “{C}entral limit theorem for {$U$}-statistics of
  associated random variables” [{S}tatist. {P}robab. {L}ett. 57 (1) (2002)
  9–15].
\newblock {\em Statist. Probab. Lett.\/}, 106:147 -- 148.

\bibitem[{Esary et~al.(1967)Esary, Proschan, and Walkup}]{MR0217826}
Esary J.D., Proschan F., and Walkup D.W. (1967).
\newblock Association of random variables, with applications.
\newblock {\em Ann. Math. Statist.\/}, 38(5):1466--1474.

\bibitem[{{Garg} and {Dewan}(2015)}]{Garg2015209}
{Garg} M. and {Dewan} I. (2015).
\newblock On asymptotic behavior of {$U$}-statistics based on associated random
  variables.
\newblock {\em Statist. Probab. Lett.\/}, 105:209 -- 220.

\bibitem[{Lebowitz(1972)}]{lebowitz1972}
Lebowitz J.L. (1972).
\newblock Bounds on the correlations and analyticity properties of
  ferromagnetic ising spin systems.
\newblock {\em Comm. Math. Phys.\/}, 28(4):313--321.

\bibitem[{Lee(1990)}]{1990u}
Lee A. (1990).
\newblock {\em U-Statistics: Theory and Practice\/}.
\newblock Statistics: A Series of Textbooks and Monographs. Taylor \& Francis.

\bibitem[{Newman(1980)}]{newman1980}
Newman C.M. (1980).
\newblock Normal fluctuations and the fkg inequalities.
\newblock {\em Comm. Math. Phys.\/}, 74(2):119--128.

\bibitem[{Newman(1984)}]{MR789244}
Newman C.M. (1984).
\newblock Asymptotic independence and limit theorems for positively and
  negatively dependent random variables.
\newblock In {\em Inequalities in statistics and probability ({L}incoln,
  {N}eb., 1982)\/}, volume~5 of {\em IMS Lecture Notes Monogr. Ser.\/}, pages
  127--140. Inst. Math. Statist., Hayward, CA.

\bibitem[{Oliveira(2012)}]{oliveira2012asymptotics}
Oliveira P. (2012).
\newblock {\em Asymptotics for {A}ssociated {R}andom {V}ariables\/}.
\newblock Springer.

\bibitem[{Prakasa~Rao(2012)}]{MR3025761}
Prakasa~Rao B.L.S. (2012).
\newblock {\em Associated sequences, demimartingales and nonparametric
  inference\/}.
\newblock Birkh\"auser/Springer, Basel.

\bibitem[{Shao and Yu(1996)}]{shao1996}
Shao Q.M. and Yu H. (1996).
\newblock Weak convergence for weighted empirical processes of dependent
  sequences.
\newblock {\em Ann. Probab.\/}, 24(4):2098--2127.

\end{thebibliography}

\end{document}